\documentclass[12pt]{article}

\usepackage[utf8]{inputenc}
\usepackage[english]{babel}
\usepackage{amsmath,amssymb,amsthm}
\usepackage{graphicx}

\def\bbC{{\mathbb{C}}}

\def\bfx{{\mathbf{x}}}
\def\calM{{\mathcal{M}}}
\def\calO{{\mathcal{O}}}
\def\calP{{\mathcal{P}}}
\def\calQ{{\mathcal{Q}}}
\def\rmd{{\mathrm{d}}}
\def\rmp{{\mathrm{p}}}
\def\dt{{\Delta t}}
\def\ds{\displaystyle}
\def\diag{\operatorname{diag}}

\newtheorem{Remark}{Remark}
\newtheorem{Definition}{Definition}
\newtheorem{Theorem}{Theorem}

\title{Strang splitting schemes for $N$-level Bloch models}
\author{Marc E. Songolo\textsuperscript{a}$^{\ast}$\thanks{$^\ast$Corresponding author. Email: marc.songolo@gmail.com} and 
Brigitte Bid\'egaray-Fesquet\textsuperscript{b} \\
\small\textsuperscript{a}D\'epartement de Math\'ematiques et Informatique, \\
\small Universit\'e de Lubumbashi, Lubumbashi, RD Congo; \\
\small \textsuperscript{b}Univ. Grenoble Alpes, Grenoble INP$^{\dag}$\footnote{$^{\dag}$Institute of Engineering Univ. Grenoble Alpes}, LJK, 38000 Grenoble, France}
\date{August 3, 2018}

\begin{document}

\maketitle

{\sc Abstract.}
We extend to the $N$-level Bloch model the splitting scheme which use exact numerical solutions of sub-equations. These exact solutions involve matrix exponentials which we want to avoid to calculate at each time step. We use Newton interpolation to reduce the computational cost. The resulting scheme is nonstandard and preserves all qualitative properties of the Bloch equations. We show numerical simulations to compare this approach with a few other schemes. 


\section{Introduction}

The spectacular development of new sources of electromagnetic radiation, such as the laser, has greatly renewed interest in studying radiation--matter interactions \cite{CDG96}. Because the intensities and the pulse times accessible by the lasers make it possible to reach the order of magnitude of the energy of cohesion of the electrons in the atoms. Some phenomena such as spontaneous or stimulated Raman scattering, Brillouin or Rayleigh, laser effect, two-photon absorption, etc. require a semi-classical model, i.e., can only be modeled with a classical field coupled with a quantum medium. Such a model is more precise than classical optics models and allows to obtain information on the structure of the atoms, thus highlighting the whole phenomenon \cite{BF06}.

In this context, the Maxwell--Bloch equations are used. The electromagnetic field is modeled by Maxwell equations and the matter is described in quantum mechanics by the Bloch equations. The coupling between these two systems is done by the expression of the polarization. The Bloch equations describe the evolution of the density matrix. They are derived from the Schrödinger equation or in the Heisenberg formalism. The density matrix is a quantum observable (unlike the wave function) and allows to describe the probability of the presence of electrons in the quantified energy levels (diagonal elements of the matrix) and the coherences between these levels (off--diagonal elements). Its size depends on the number of levels. In many references, the derivation of the Bloch equation is only presented in the case of two-level atoms.

Bloch equations can already raise problems both from a theoretical and a numerical point of view. Indeed, they have been approached by different methods such as the Crank--Nicolson method \cite{ ZAG95, Zio97a, Zio97b, NY98, GS03, AC09, Bid03}, a fourth-order Runge--Kutta method \cite{Uwi15}, the relaxation method \cite{CN04}, etc. But most of these numerical schemes do not conserve the qualitative properties of Bloch equations. To overcome this deficiency, Bidégaray \textit{et al.} \cite{BBR01} have introduced a Strang splitting method which preserve some physical properties (Hermicity, positivity, trace conservation).

Following these study, we have improved Strang splitting schemes for the two-level Bloch model in \cite{SBF17}. A special feature of these schemes is that the solutions of sub-equations are exact, of variable time step-size for the Liouville equation and conform to the nonstandard finite difference (NSFD) methodology developed by Mickens \cite{Mic94}. These splitting schemes conserve all the physical properties of the Bloch equations. Moreover, they are explicit and retain the advantage of stability when coupled with Maxwell equations. In this paper, we want to extend this type of schemes to the $ N $-level Bloch model, hoping thus to improve the scheme presented in \cite{BBR01}.

The paper is organized as follows: we introduce the Bloch model in Section 2. Section 3 provides the construction rules of NSFD schemes. In Section 4, we first present the decomposition of the Bloch equation, then the exact schemes of the sub-equations, and finally the splitting scheme. In order to reduce the algorithmic complexity of the splitting scheme, we propose in Section 5 an equivalent formulation of the matrix exponentials that occur in the exact discretization of the Liouville equation. In Section 6, we compare numerically the Crank--Nicolson method and the reformulations of the exponential discussed in this paper. The study is followed by an appendix containing an alternative method for the derivation of the matrix exponential.

\section{Bloch model}

The derivation of the Bloch equation can be found in many textbooks (see for example \cite{BF06, Blu12, Boy08, Lou00, WH03}). In this study, we use dimensionless equations:
\begin{equation}
\label{eq:BlochSansDim}
\partial_t\rho = -i [H_0,\rho] - i [V,\rho] + Q(\rho),
\end{equation}
where $[A,B]$ is the commutator of the two operators $A$ and $B$. The diagonal elements of the \textit{density matrix} $\rho$ are called \textit{populations} and express the presence probabilities of electrons in the quantified energy levels. The off-diagonal elements are complex numbers called \textit{coherences}, whose moduli can be interpreted as conditional probabilities of transition between the energy levels. 
In equation \eqref{eq:BlochSansDim}, $H_0$ is the free Hamiltonian of electrons and is a diagonal matrix $\diag(\hbar\omega_j)_{j=1,\dots,n}$.
The potential $V(t)$ is a zero diagonal, Hermitian matrix and results from the interaction with an electromagnetic wave.  
A phenomenological relaxation term $Q(\rho)_{jk}$ can be added to model many phenomena such as spontaneous emission, collisions, vibrations in crystal lattices, etc. It must be chosen so as to preserve over time some properties of the density matrix, in particular Hermicity, positiveness and trace conservation (see \cite{BBR01} for details).

\section{NSFD Schemes}

To discuss NSFD schemes, we cast the Bloch equation as
\begin{equation}
\label{eq:formred}
\partial_{t}\rho = F(\rho),
\end{equation}
with unknown the matrix $\rho:[t_0,T]\rightarrow \bbC^{n\times n}$, initially equal to $\rho_0\in\bbC^{n\times n}$, and $F:\bbC^{n\times n}\rightarrow\bbC^{n\times n}$ is a given function.

For the numerical approximation of \eqref{eq:formred}, we discretize the interval $[t_0,T]$ at the discrete times $t_n=t_0+n\dt$, where the parameter $\dt>0$ is the step-size. We denote by $\rho^n$ an approximation of the solution $\rho(t_n)$ at time $t_n$.

The finite difference equation reads
\begin{equation}
\label{eq:formNSFDS}
D_\dt(\rho^n) = F_\dt(t_n,\rho^n),
\end{equation}
where $D_\dt(\rho^n)$ and $F_\dt(t_n,\rho^n)$ approximate $\partial_t\rho(t_n)$ and $F(t_n,\rho(t_n))$ respectively.

\begin{Definition}
The scheme \eqref{eq:formNSFDS} is called a nonstandard finite difference method if at least one of the following conditions is satisfied:
\begin{itemize}
\item $D_\dt(\rho^n) = (\phi(\dt))^{-1}(\rho^{n+1}-\rho^n)$ where $\phi(\dt) = \dt I +\calO(\dt^2)$ is a positive diagonal matrix;
\item $F_{\dt}(t_n, \rho^n)=g(\rho^{n},\rho^{n+1},\dt)$ where $g(\rho^n,\rho^{n+1},\dt)$ is a nonlocal approximation of the right-hand side of System \eqref{eq:formred}.
\end{itemize}
\end{Definition}

These notions are discussed in detail in \cite{Mic94, Mic00, Mic05, AL03}. Moreover, Mickens has introduced in \cite{Mic00} a rule for the construction of NSFD schemes for complex equations.\\
\textbf{\textit{Rule for complex equations.}}
\textit{For differential equations having $N$ ($\geq 3$) terms, it is generally useful to construct finite difference schemes for various sub-equations composed of $M$ terms, where $M<N$, and then combine all the schemes together in an overall consistent finite difference model.}

By this last rule, it is necessary to split Equation \eqref{eq:BlochSansDim} into two sub-equations, then solve sub-equations by exact methods, and finally, connect solutions of sub-equations through a single consistent solution. To this end, we explore how to construct consistent finite difference models using Strang splitting method.

\section{Splitting method}

We rewrite Equation \eqref{eq:BlochSansDim} as
\begin{equation*}
\partial_t\rho_{jk} = -i\omega_{jk}\rho_{jk} - i[V,\rho]_{jk} + Q(\rho)_{jk},
\end{equation*}
where $\omega_{jk} = \omega_j - \omega_k$ is the frequency associated with the transition from level $k$ to level $j$. This equation is decomposed into the relaxation--nutation evolution
\begin{equation}
\label{eq:EqRelaxNut}
\partial_t\rho = L\rho,
\end{equation}
where $(L\rho)_{jk} = -i\omega_{jk}\rho_{jk} + Q(\rho)_{jk}$ and the interaction with the electromagnetic field
\begin{equation}
\label{eq:Liouville}
\partial_t\rho_{jk} = - i[V,\rho]_{jk}.
\end{equation}
We have seen in \cite{SBF17} that this splitting yields the best approximation for the Bloch equation, and even is better than no splitting for a Self-Induced Transparency test case.

As the relaxation--nutation operator is linear and time invariant, the solution of Equation \eqref{eq:EqRelaxNut} is
\begin{equation}
\label{eq:SolRelaxNut}
\rho(t) = \exp(L(t-t_0))\rho(t_0).
\end{equation}
Since the potential $V$ generally depends on time, the solution of the interaction equation \eqref{eq:Liouville} is
\begin{equation}
\label{eq:SolLiouville}
\rho(t) = \exp\left(- i\int_{t_0}^t V(\tau)\rmd\tau\right)\rho(t_0)
\exp\left(i\int_{t_0}^t V(\tau)\rmd\tau\right).
\end{equation}

\subsection{Exact discretization of the relaxation--nutation equation}

An exact finite difference scheme for Equation \eqref{eq:EqRelaxNut} is easily deduced from its analytical solution \eqref{eq:SolRelaxNut} and one time-step of the relaxation--nutation equation reads
\begin{equation}
\label{eq:SchemaRelaxNut}
\rho^{n+1} = e^{L\dt} \rho^n.
\end{equation}

\subsection{Exact discretization of the Liouville equation}

Let $V^{n+1/2}$ be the mean of the function $V$ on the interval $[t_n,t_{n+1}]$:
\begin{equation*}
V^{n+1/2} = \frac1\dt \int_{t_{n}}^{t_{n+1}}V(\tau)\rmd\tau.
\end{equation*}
Then one time-step of the Liouville equation is easily deduced from 
\eqref{eq:SolLiouville}, namely
\begin{equation}
\label{eq:SchemaLiouville}
\rho^{n+1} = \exp(-i\dt V^{n+1/2}) \rho^n \exp(i\dt V^{n+1/2}).
\end{equation}

\subsection{Strang splitting}

To construct the splitting scheme, we choose the Strang formula \cite{Str68} in order to achieve second order precision, which would prove useful when coupling with an order-two scheme for the electromagnetic field in a Maxwell--Bloch context (see \cite{BF06}).

Furthermore, this method is consistent (see \cite{FH07}, for details) according to the discretization rule for complex equations, and preserve the symmetry and positiveness properties of the density matrix \cite{BF06}, provided relaxation terms satisfy the conditions given in \cite{BBR01}.
Thus, the Strang splitting method for the Bloch model reads
\begin{equation}
\label{eq:Schemasplitting}
\rho^{n+1} 
= \exp(L\dt/2)\exp(-i\dt V^{n+1/2}) \exp(L\dt/2)\rho^n \exp(i\dt V^{n+1/2}).
\end{equation}
Here, we begin and end up with the relaxation--nutation term in the splitting. This is the good choice, since this term is the steepest one when the electromagnetic field is small, which necessarily happens in the test cases
(oscillating field, wave paquet).
In the sequel, we give equivalent formulations to $\exp(i\dt V^{n+1/2})$, to avoid the possible complexity of the calculation of matrix exponentials at each time step.

\section{Exponential of $N\times N$ matrices}
\label{sec:matricesNN}

Dozens of methods for calculating the exponential of a matrix can be obtained from the more or less classical results in analysis, approximation theory, and matrix theory. In \cite{ML03}, the authors describe nineteen methods that seem to be practical. The relative effectiveness of each method is evaluated according to the following attributes, listed in decreasing order of importance: generality, reliability, stability, accuracy, efficiency, storage requirements, ease of use and simplicity. Generality means that the process is applicable to large classes of matrices. For example, a method that only works on matrices with distinct eigenvalues will not be much appreciated. Moreover, an algorithm will be said to be stable if it has no more sensitivity to perturbations than is inherent to the underlying problem. The precision of an algorithm refers mainly to the error introduced by truncating an infinite series to a certain order. By these standards, none of the algorithms we know is satisfactory, although some are much better than others. 

\subsection{Interpretation as an interpolation problem}

The Cayley--Hamilton theorem applied to a matrix $A\in\calM_{N}(\bbC)$ ensures that $p(A)=0$ where $p$ is the characteristic polynomial of $A$ defined by $p(\lambda)=\det(\lambda I-A)$.
This allow to express $A^N$ and higher powers as $N-1$-degree polynomials of $A$. More generally it allows to express an analytical function of $A$ as such a polynomial. This is the case for the exponential of $A$. \\
Let $\gamma\in\bbC$, the function $\exp(i\gamma A)$ can also be expressed by a polynomial that we denote $\calP_\gamma(A)$:
\begin{equation}
\label{eq:expNewton}
\calP_\gamma(A)=\sum_{j=0}^{N-1}\alpha_{j}(\gamma) A^j.
\end{equation}
If the eigenvalues $\lambda_1,\dots,\lambda_N$ of matrix $A$ are distinct, then there is a basis in which $A$ is the diagonal matrix $D$. Let $P$ be the change of basis matrix, which does not depend on $\gamma$,  we have
\begin{equation*}
\exp(i\gamma D) = P^{-1}\exp(i\gamma A)P = P^{-1}\calP_\gamma(A)P 
= \calP_\gamma(D).
\end{equation*}
This relation involving diagonal matrices is actually a system of $N$ independent equations 
\begin{equation*}
\exp(i\gamma \lambda_k) = \calP_\gamma(\lambda_k),\ k=1,\dots,N,
\end{equation*}
which admits a unique solution, because it amounts to inverting a Vandermonde matrix. We however do not want to invert it, but interpret this system as an interpolation problem, i.e. interpolate function $\calP_\gamma$ at the locations $\lambda_k$, $k=1,\dots,N$, with values $\exp(i\gamma\lambda_k)$. Since the function $\calP_\gamma$ is an $(N-1)$-degree polynomial, the interpolation polynomial at these $N$ locations will be exactly the function itself.

\subsection{Newton interpolation}

In \cite{SBF17}, in the case of 2-level Bloch equations, interpolation polynomial is expressed in the canonical basis $(I,A)$. This can of course be extended to the case of $N\times N$ matrices in the $(I,A,\dots,A^{N-1})$ basis but the formulae are quite intricate and have to be derived individually for each value of $N$. The $3\times3$ case is treated in an Appendix.\\
 
Therefore we use here the Newton basis
\begin{equation*}
I, (\gamma A-\lambda_1I), (\gamma A-\lambda_1 I)(\gamma A-\lambda_2I), \dots, (\gamma A-\lambda_1 I)\dots(\gamma A-\lambda_{N-1}I).
\end{equation*}
In this decomposition the coefficients are divided differences
\begin{equation*}
\calP_\gamma(A) = f[\lambda_1] + \sum_{\ell=2}^{N} f[\lambda_1,\dots,\lambda_\ell] \prod_{k=1}^{\ell-1}(\gamma A-\lambda_k I),
\end{equation*}
where $f[\lambda_k] = \exp(i\gamma\lambda_k)$, $k=1,\dots,N$ and we have the recursion formula
\begin{equation*}
f[\lambda_k,\dots,\lambda_\ell] = \frac{f[\lambda_k,\dots,\lambda_{\ell-1}] - f[\lambda_{k+1},\dots,\lambda_\ell]}{\lambda_k-\lambda_\ell},\ 1\leq k<\ell\leq N.
\end{equation*}
The calculation of the polynomial is then done iteratively using the Hörner algorithm. Indeed, we have
\begin{align}
\calP_\gamma(A) & = c_0 + \sum_{\ell=1}^{N-1} c_\ell \prod_{k=1}^{\ell}(\gamma A-\lambda_k I) \nonumber \\
\label{eq:Horner}
& = c_0 + (\gamma A-\lambda_1)(c_1 + (\gamma A-\lambda_2)(c_2+\dots+(\gamma A-\lambda_{N-1}))).
\end{align}
The advantage of this approach, compared to the one based on the canonical basis, is that the numerical code produced is generic for all $N$. 

\subsection{NSFD interpretation for the Liouville equation}

In Equations \eqref{eq:SchemaLiouville} and \eqref{eq:Schemasplitting} we need to evaluate $\exp(i\dt V^{n+1/2})$. In the dimensionless version of the equations we are dealing with, matrix $V^{n+1/2}$ is the product of the scalar electric field $E^{n+1/2}$ and a constant polarisability matrix $p$ (see \cite{BF06}). We can therefore write
\begin{equation*}
\exp(i\dt V^{n+1/2})=\calP_{\dt E^{n+1/2}}(p).
\end{equation*}
Therefore, the exact scheme for the Liouville equation can be written as
\begin{equation}
\label{eq:exactLiouville}
\rho^{n+1} = \calP_{\dt E^{n+1/2}}^{-1}(p) \rho^n \calP_{\dt E^{n+1/2}}(p),
\end{equation}
or equivalenty
\begin{equation}
\label{eq:exactLiouville2}
\calP_{\dt E^{n+1/2}}(p) \rho^{n+1} = \rho^n \calP_{\dt E^{n+1/2}}(p).
\end{equation}
The  polynomial $\calP_{\dt E^{n+1/2}}(p)$ is of course also equal to the series expansion of the exponential $\exp(i\dt V^{n+1/2})$. Therefore, in the limit $\dt\to 0$, $\alpha_j(\dt E^{n+1/2}) = (i\dt E^{n+1/2})^j/(j!) + O(\dt^{j+1}) $. In particular $\alpha_0(\dt E^{n+1/2}) = 1 + O(\dt)$ and $\alpha_1(\dt E^{n+1/2}) = i\dt E^{n+1/2} + O(\dt^2)$. Let us set
$\alpha_1(\dt E^{n+1/2}) = iE^{n+1/2}\tilde\alpha_1(\dt E^{n+1/2})$. For a small enough $\dt$, we can therefore ensure that $\alpha_0(\dt E^{n+1/2})$ and $\tilde\alpha_1(\dt E^{n+1/2})$ are non zero and rewrite Equation \eqref{eq:expNewton} as
\begin{equation*}
\calP_{\dt E^{n+1/2}}(p) = \alpha_0(\dt E^{n+1/2}) + \tilde\alpha_1(\dt E^{n+1/2}) i E^{n+1/2} \calQ_{\dt E^{n+1/2}}(p), 
\end{equation*}
where $\calQ_{\dt E^{n+1/2}}(p)$ is an order $N-1$ polynomial with smaller degree term equal to $p$. The exact scheme for the $N$-level Liouville equation \eqref{eq:exactLiouville2} reads
\begin{equation}
\label{eq:schemaimpliciteN}
\left\{
\begin{array}{l}
\begin{aligned}
(\Phi^{n+1/2}(\dt))^{-1}(\rho^{n+1}-\rho^{n}) = & -i \left\{\tilde V^{n+1/2}\rho^{n+1}-\rho^{n}\tilde V^{n+1/2}\right\},\\
\Phi^{n+1/2}(\dt) = & \dfrac{\tilde\alpha_1(\dt E^{n+1/2})}{\alpha_0(\dt E^{n+1/2})}I,
\end{aligned}
\end{array}
\right.
\end{equation}
where
\begin{equation*}
\tilde V^{n+1/2} = E^{n+1/2} \calQ_{\dt E^{n+1/2}}(p).
\end{equation*}

\begin{Remark}
Observe a nonlocal discretization and a renormalisation of the step-size in the exact scheme for the $ N $-level Liouville equation, according to the Mickens rules. In particular $\Phi(\dt) = \dt I +\calO(\dt^2)$.
\end{Remark}
\begin{Remark}
$\tilde V^{n+1/2}$ is not just equal to $V^{n+1/2}$ but there is an additional higher order term which is called a recovery factor because it comes from estimating the matrix exponential. The addition of this term contradicts one of the basic principles of the NSFD schemes theory, which forbids any form of adjustment by adding ad hoc terms. 
\end{Remark}

\section{Numerical simulations}

The decomposition of the Bloch equation into the relaxation--nutation evolution and the evolution of interaction with the electromagnetic field has two main advantages. First, part of the computations can be performed off-line, before the time iterations, calculating once and for all the eigenvalues of the electric dipole matrix, as well as the change of basis matrix and its inverse. Here we deal only with the Bloch equation, but this is even more efficient when there is space dependence, such as when coupling with Maxwell equations. The second advantage has been pointed out in \cite{SBF17} and demonstrated on a Self Induced Transparency test case, and is the fact that it can decouple stiff and non-stiff parts of the equation.

In the following, we compare various schemes for the Bloch model, always performing the same decomposition but varying the way the exponential is calculated or approximated. We consider the historical method for the Bloch equation \cite{ZAG95,Zio97a,Zio97b}, namely the Crank--Nicolson method, although it has been shown in Reignier's thesis \cite{Rei00} that it does not preserve the property of positivity for more than three levels. Other methods such as the fourth-order Runge--Kutta method are not adapted to preserve the physical properties of the Bloch equation \cite{SBF17}. We compare the Crank--Nicolson method with the method presented in this paper, and the computation of the matrix exponential.
We describe below the four methods.  

\paragraph*{\textit{Exponential method.}} 

\begin{equation*}
\rho^{n+1} = \exp(L\dt/2) \exp(-i\dt V^{n+1/2}) \exp(L\dt/2) \rho^n \exp(i\dt V^{n+1/2}).
\end{equation*}

\paragraph*{\textit{Crank--Nicolson method.}} 

The matrix exponential is approximated by the Crank--Nicolson scheme:
\begin{equation}
\left\{
\begin{array}{l}
\rho^{n+1} = \exp(L\dt/2) (A^{n+1/2})^{-1} \exp(L\dt/2) \rho^n A^{n+1/2},\\
\ds A^{n+1/2} = \left(I+\frac{i}2\dt V^{n+1/2}\right)\left(I-\frac{i}2\dt V^{n+1/2}\right)^{-1}.
\end{array}
\right.
\end{equation}

\paragraph*{\textit{Newton method.}} 

The method described in this paper is used, the Newton basis is used and the polynomial is reconstructed \textit{via} the Horner algorithm:
\begin{equation}
\left\{
\begin{array}{l}
\rho^{n+1} = \exp(L\dt/2) (B^{n+1/2})^{-1} \exp(L\dt/2) \rho^n B^{n+1/2},\\
\ds B^{n+1/2} = \calP_{\dt E^{n+1/2}}(p).
\end{array}
\right.
\end{equation}

\paragraph*{\textit{Canonical method.}} 

As in \cite{SBF17} the polynomial equivalent to the exponential is expressed in the canonical form. Details for three levels can be found in Appendix \ref{app:canonique}:
\begin{equation}
\label{eq:canonicalscheme}
\left\{
\begin{array}{l}
\rho^{n+1} = \exp(L\dt/2) (C^{n+1/2})^{-1} \exp(L\dt/2) \rho^n C^{n+1/2},\\
\ds C^{n+1/2} = \sum_{j=0}^{N-1} \alpha_{j}^{n+1/2} (V^{n+1/2})^{j}.
\end{array}
\right.
\end{equation}

\subsection{Three-level test case}

\begin{figure}[h]
\centering
\includegraphics[width=.55\textwidth]{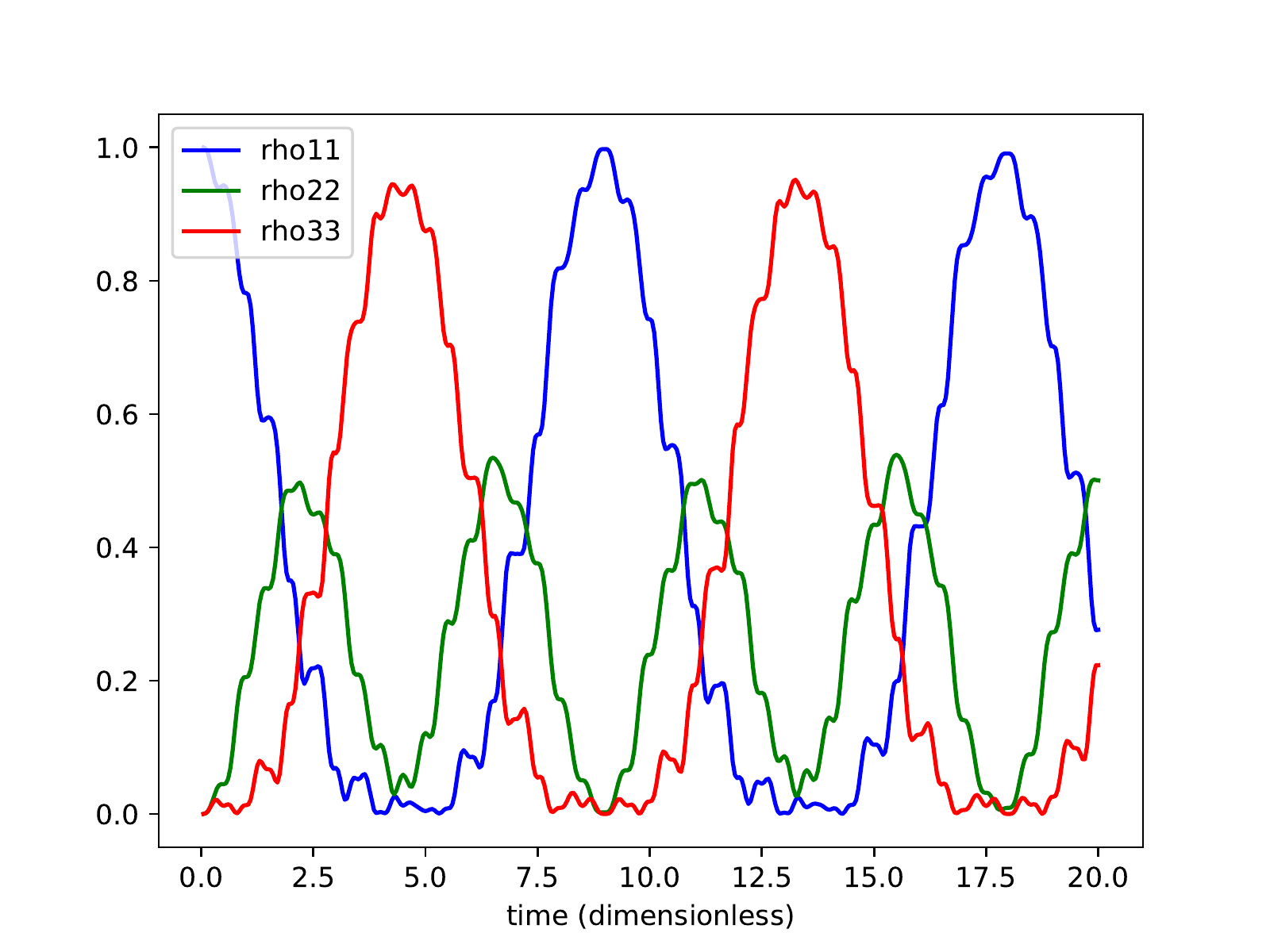}
\caption{\label{fig:test3N} Time evolution of populations over 20 periods for a 3-level test-case.}
\end{figure}

We first compare the methods on a three-level test case. We suppose there is no relaxation ($Q(\rho)=0$) and apply a sinusoidal input electrical field $E(t)=\sin(2\pi t)$ (recall we deal with dimensionless equations). The level frequencies are chosen to be $0$, $\pi$ and $2\pi$ since resonance with the input wave is required for the system to evolve nontrivially. Besides the polarizability matrix $p$ is chosen to be
\begin{equation*}
p = \begin{pmatrix}
0 & 1 & 1.1 \\
1 & 0 & 1 \\
1.1 & 1 & 0
\end{pmatrix}.
\end{equation*}
Let $n_\rmp$ be the number of discretization times within one period of the input signal. The time step is therefore equal to $1/n_\rmp$. The time evolution of populations over 20 periods of the input signal is displayed on Figure \ref{fig:test3N}. This result has been obtained with the Python implementation of the matrix exponential $\exp(i\dt V^{n+1/2})$, and $n_\rmp=20$, but similar results can be obtained with the other methods described in this paper.

Notice that the pseudo-period of the population is not that of the wave. The period of the wave can be however seen on the plots, since every flat parts correspond to a vanishing input electric field. 

Now various parameters will vary, to begin with the numerical method and the time step. We compare the computational time (for a not especially optimized python implementation on a small laptop). The results are gathered in Table \ref{tab:test3N}. To homogenize the results, 2000 periods have been simulated (which, coming back to a dimensional world, would correspond for light waves to about 10 ps). 

\begin{table}[h]
\centering
\begin{tabular}{c|cccc}
$n_\rmp$ & Exponential & Crank-Nicolson & Newton & Canonical \\
\hline
5 & 10 & 4 (out) & 4 & 3 \\
10 & 14 & 6 (bad) & 7 & 5 \\
20 & 27 & 11 (bad) & 15 & 11 \\
100 & 131 & 38 & 78 & 53
\end{tabular}
\caption{\label{tab:test3N} Computational time (in seconds) for a 3-level test-case.}
\end{table}

The Exponential method is clearly the most expensive. The Crank-Nicolson is the least expensive, but the quality of the results disqualifies this method since we need a lot of points to ensure the same quality as with the other methods. With $n_\rmp=5$, positiveness is violated from the very first periods on. Therefore, the Newton and Canonical methods seem to be the best, from the computational time point of view, with a little advance for the Canonical methods. We see next why we however prefer the Newton method.

\begin{figure}
\begin{center}
\begin{tabular}{c}
\includegraphics[width=.55\textwidth]{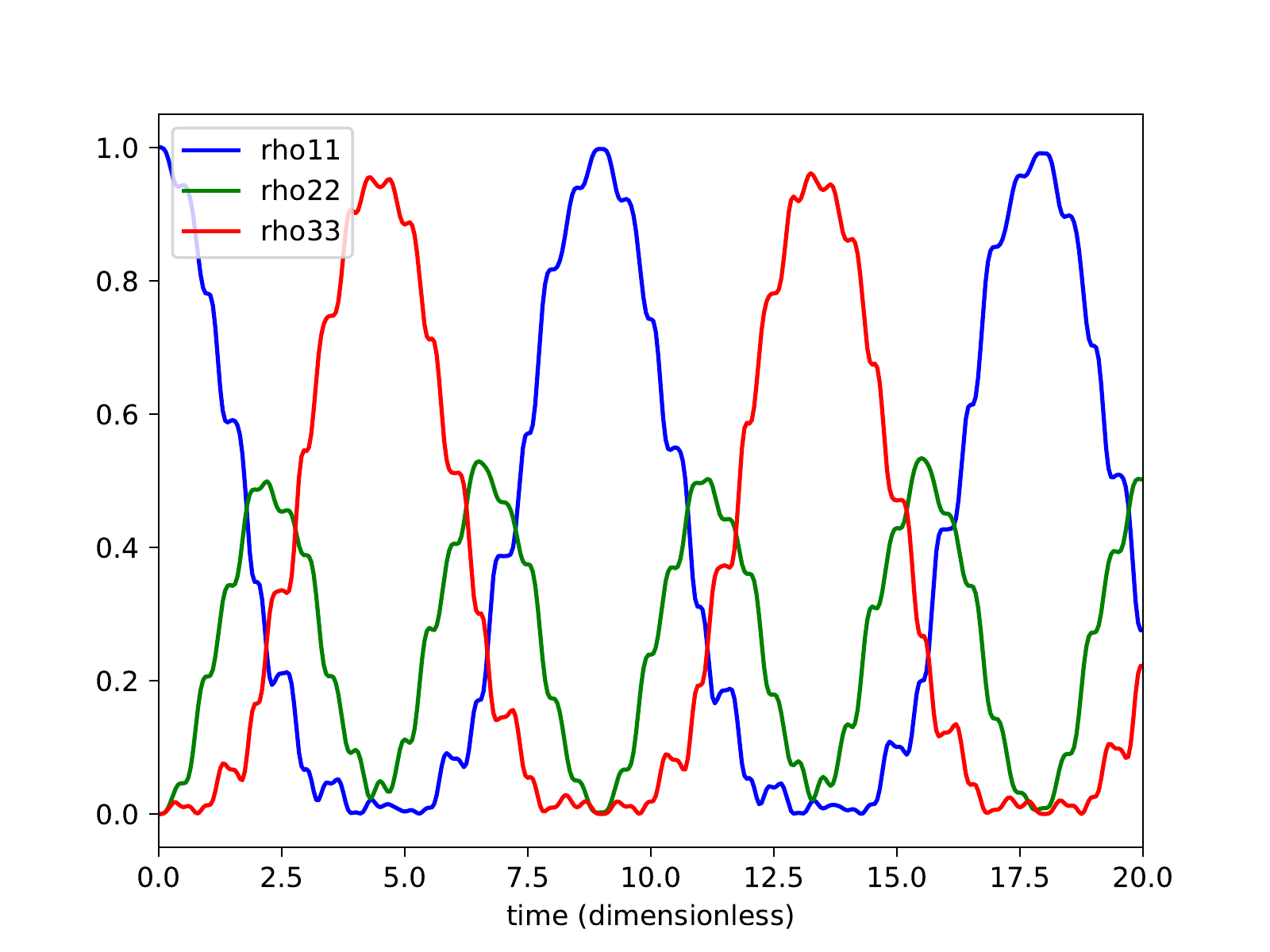}\\
\includegraphics[width=.55\textwidth]{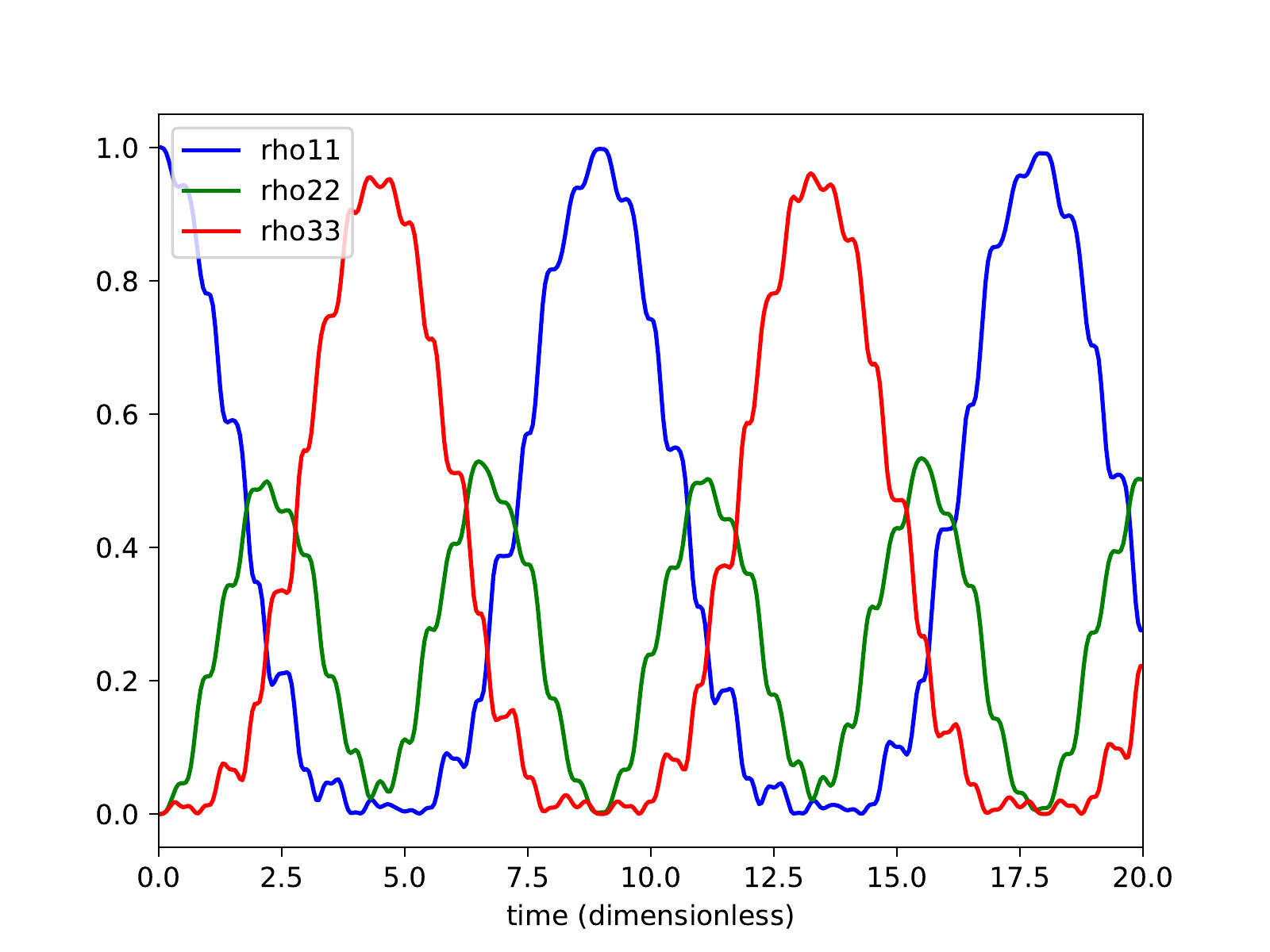}\\
\includegraphics[width=.55\textwidth]{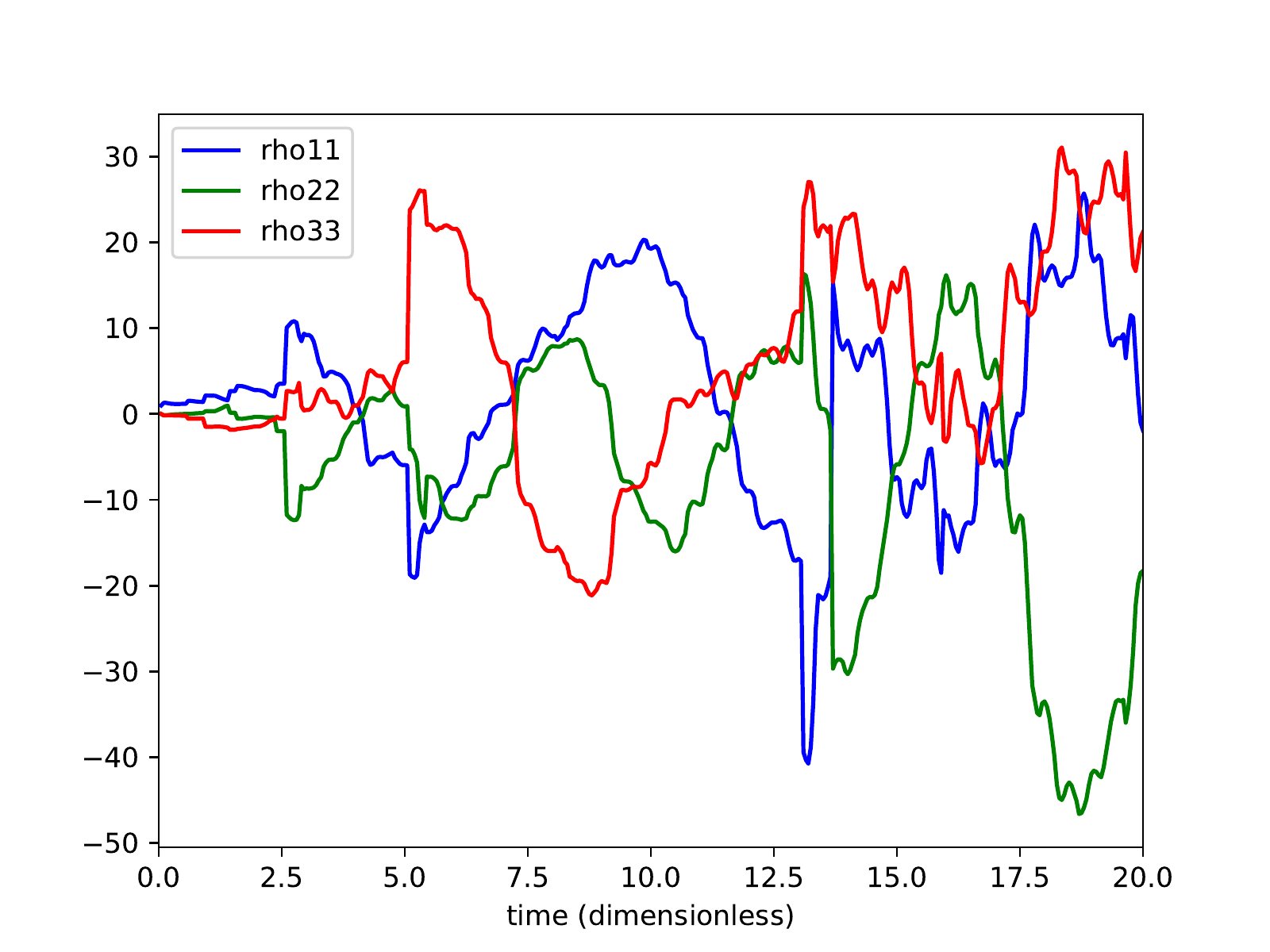}\\
\end{tabular}
\caption{\label{fig:test3N_degenerate} Time evolution of populations over 20 periods in a degenerate situation. From top to bottom: Exponential, Newton, and Canonical method.}
\end{center}
\end{figure}

We have previously chosen a strange matrix $p$ to prevent it to have equal eigenvalues. Let us now take \begin{equation*}
p = \begin{pmatrix}
0 & 1 & 1 \\
1 & 0 & 1 \\
1 & 1 & 0
\end{pmatrix},
\end{equation*}
which eigenvalues are $-1$ (double) and 2, and come back to 20 periods of the input signal and $n_\rmp=20$. The results are displayed on Figure \ref{fig:test3N_degenerate} for the Newton and the Canonical methods. The Exponential method serves here as a reference solution. 

While both the Newton and the Canonical methods use the eigenvalues of the problem to compute the matrix exponential, the Canonical method is very sensitive to degenerate situations. In this case, the Jordan form of the the matrix is not diagonal and different formulae should be used (see \cite{DH18} for full details). This does not solve the problem, since in the case of a nearly degenerate situation (two very close eigenvalues) the formulae presented in this paper (see Appendix) should used but would be very unstable.

\subsection{$N$-level test case}

Now, we compare only the Exponential and Newton methods, and have $N$ vary. The results, namely the computational times, are gathered in Table \ref{tab:testNN} for various values of $N$ and anew 2000 periods of the input signal.

\begin{table}[h]
\centering
\begin{tabular}{c|cccc}
$N$ & Exponential & Newton \\
\hline
2 & 27 & 9 \\
3 & 27 & 16 \\
4 & 27 & 23\\
5 & 27 & 34 \\
10 & 34 & 125 
\end{tabular}
\caption{\label{tab:testNN} Computational time (in seconds) for a N-level test-case.}
\end{table}

We are not very fair with the Newton method, since we compare it with a clearly well optimized Python matrix exponential which computational time barely depends on the dimension of the matrix. In many studies however we are dealing with small density matrices, describing 2, 3 or 4 levels. In these cases the computational gain using the Newton method is not impressive, but can prove very useful for simulation over long physical times or involving many space locations.

\section{Conclusion}

We have derived splitting schemes for the $N$-level Bloch model. To this aim, the Bloch equation has been decomposed into a relaxation--nutation evolution and the interaction with the electromagnetic field (which is a Liouville equation). We are able to obtain exact solutions for the resulting sub-equations, and construct a Strang splitting scheme. The solution of the Liouville equation involves matrix exponentials and we discuss whether it is reasonable or not to compute it. Indeed, thanks to the Cayley--Hamilton theorem, it can be replaced by the computation of a polynomial. We used in particular Newton interpolation to define this polynomial. The resulting scheme has a variable time step-size and satisfies the nonstandard discretization rules of Mickens. Moreover, the splitting scheme preserve the qualitative and quantitative properties (Hermicity, trace conservation, positiveness) of Bloch equations.

The numerical comparison of the methods show that computing a polynomial instead of the exponential is advantageous for small density matrices, i.e. a small number of levels. If the gain is relatively low, the number of such computations for a full Maxwell--Bloch simulation can really make it a good track to improve the computational load.  

The choice of the splitting scheme aims at solving exactly each sub-equation, but also at dealing correctly with terms with different stiffness. An interesting goal in this direction is to derive methods that preserve the asymptotic behavior to the rate equations. If splitting methods are not direct candidates for this since they dissociate different parts of the Bloch equations which are intimately connected in the Boltzmann equation, it would be interesting to connect numerical solutions of these two sub-equations in an implicit method by the NSFD technique, hoping thus to obtain asymptotic preserving schemes. This is the object of our future research.

\appendix

\section{Alternative method for $3\times3$ matrix exponentials}
\label{app:canonique}
\setcounter{equation}{0}
\def\theequation{A.\arabic{equation}}

\subsection{Problem setting}

The problem of the expression of exponential of matrices as exact finite difference schemes has been studied in \cite{DH18} for general $3$-order systems
\begin{equation}
\label{eq:systlineairetrois}
\bfx'(t) = A\bfx(t);\ \bfx(t)=(x_1(t),x_2(t),x_3(t))^T,\ A\in\calM_3(\bbC).
\end{equation}
Here we restrict to the case where matrix $A$ is similar to the canonical Jordan form
\begin{equation}
\label{eq:formejordan3}
J_3 = \begin{pmatrix}\lambda_1 & 0 & 0\\ 0 & \lambda_2 & 0\\0 & 0 & \lambda_3 \end{pmatrix},
\end{equation}
where $\lambda_1$, $\lambda_2$, and $ \lambda_3$ are distinct.\\
We have already seen that solving numerically $\bfx'(t) = J_3\bfx(t)$ amounts to solving three decoupled systems
\begin{equation}
\label{eq:exactnumjordantrois}
x_1^{n+1} = x_1^n e^{\lambda_1\dt},\ 
x_2^{n+1} = x_2^n e^{\lambda_2\dt},\
x_3^{n+1} = x_3^n e^{\lambda_3\dt}. 
\end{equation}

Back in the original basis, this can be written as
\begin{equation}
\label{eq:schemaexplicite}
\bfx^{n+1} = (\alpha(\dt) I + \beta(\dt) A + \gamma(\dt) A^2) \bfx^n.
\end{equation}
which we prefer to express as a polynomial of $\beta(\dt)A$:
\begin{equation}
\bfx^{n+1} = (\alpha(\dt) I + \beta(\dt) A + \xi(\dt)(\beta(\dt))^2 A^2) \bfx^n.
\end{equation}
The explicit exact finite difference form is
\begin{equation*}
\dfrac{\bfx^{n+1}-\alpha(\dt)\bfx^n}{\beta(\dt)}=A\bfx^n+\xi(\dt)\beta(\dt) A^{2}\bfx^n,
\end{equation*}
where $\alpha(\dt)$, $\beta(\dt)$, and $\xi(\dt) $ are parameters to be determined.
The same coefficients appear in the Jordan basis, and identifying in Equation \eqref{eq:exactnumjordantrois} we obtain
\begin{equation*}
\label{eq:Vandermonde}
\alpha(\dt) + \beta(\dt)\lambda_j + \xi(\dt)(\beta(\dt))^2\lambda_j^2 = e^{\lambda_j\dt},\ j=1,2,3.
\end{equation*}

\subsection{Coefficients}

Solving the Vandermonde system \eqref{eq:Vandermonde} is a classical problem. It can be expressed using determinants. Defining the determinant of the system:
\begin{equation*}
\delta = \begin{vmatrix}1 & \lambda_1 & \lambda_1^2 \\
1 & \lambda_2 & \lambda_2^2 \\
1 & \lambda_3 & \lambda_3^2\end{vmatrix}, 
\end{equation*}
the coefficients are equal to 
\begin{equation*}
\delta\alpha(\dt) = \begin{vmatrix}e^{\lambda_1\dt} & \lambda_1 & \lambda_1^2 \\
e^{\lambda_2\dt} & \lambda_2 & \lambda_2^2 \\
e^{\lambda_3\dt} & \lambda_3 & \lambda_3^2\end{vmatrix},
\delta\beta(\dt) = \begin{vmatrix}1 & e^{\lambda_1\dt} & \lambda_1^2 \\
1 & e^{\lambda_2\dt} & \lambda_2^2 \\
1 & e^{\lambda_3\dt} & \lambda_3^2\end{vmatrix},
\delta\gamma(\dt) = \begin{vmatrix}1 & \lambda_1 & e^{\lambda_1\dt} \\
1 & \lambda_3 & e^{\lambda_2\dt} \\
1 & \lambda_3 & e^{\lambda_3\dt} \end{vmatrix}.
\end{equation*}
More explicitely we have 
\begin{equation}
\label{eq:paraexactjordanexpl}
\begin{aligned}
\delta & = (\lambda_2-\lambda_1)(\lambda_3-\lambda_1)(\lambda_3-\lambda_2), \\
\alpha(\dt) & = \dfrac{e^{\lambda_1\dt}\lambda_2\lambda_3(\lambda_3-\lambda_2)+e^{\lambda_2\dt}\lambda_3\lambda_1(\lambda_1-\lambda_3)+e^{\lambda_3\dt}\lambda_1\lambda_2(\lambda_2-\lambda_1)}\delta,\\
\beta(\dt) & = \dfrac{e^{\lambda_1\dt}(\lambda_2^2-\lambda_3^3)+e^{\lambda_2\dt}(\lambda_3^2-\lambda_1^3)+e^{\lambda_3\dt}(\lambda_1^2-\lambda_2^3)}\delta,\\
\gamma(\dt) & = \dfrac{e^{\lambda_1\dt}(\lambda_3-\lambda_2)+e^{\lambda_2\dt}(\lambda_1-\lambda_3)+e^{\lambda_3\dt}(\lambda_2-\lambda_1)}\delta,\\
\xi(\dt) & = \frac\gamma{\beta^2}.
\end{aligned}
\end{equation}

\begin{Theorem}
\label{th:appendix}
For any matrix $A\in\calM_3(\bbC)$, whose eigenvalues are distinct,
\begin{equation}
\label{eq:expsemimpl}
\exp(\dt A) = \alpha(\dt) I + \beta(\dt) A + \xi(\dt)(\beta(\dt))^2 A^2
\end{equation}
where $\alpha(\dt)$, $\beta(\dt)$ and $\xi(\dt)$ are determined by the relations \eqref{eq:paraexactjordanexpl}.
\end{Theorem}

\subsection{Application to the three-level Bloch equation}

We now want to make explicit the coefficients $\alpha_0$, $\alpha_1$ and $\alpha_2$ in Equation
\eqref{eq:canonicalscheme} to construct matrix $C^{n+1/2}$:
\begin{equation}
\label{eq:canonicalscheme3} 
\exp(i\dt V^{n+1/2}) = \alpha_0^{n+1/2}I+\alpha_1^{n+1/2}V^{n+1/2}+\alpha_2^{n+1/2} (V^{n+1/2})^2.
\end{equation}
Compared Equation \eqref{eq:expsemimpl} of Theorem \ref{th:appendix} and Equation \eqref{eq:canonicalscheme3}, we have $A=iV^{n+1/2}$, $\alpha_0^{n+1/2}=\alpha(\dt)$, $\alpha_1^{n+1/2}=i\beta(\dt)$, $\alpha_2^{n+1/2}=-\xi(\dt)\beta^2(\dt)$. Denoting $\theta_j^{n+1/2} = \lambda_j^{n+1/2}\dt$, $j=1,2,3$, where the $\lambda_j^{n+1/2}$ are the distinct eigenvalues of matrix $V^{n+1/2}$, we obtain
\begin{equation*}
\begin{aligned}
\delta^{n+1/2} & = (\lambda_2^{n+1/2}-\lambda_1^{n+1/2})(\lambda_3^{n+1/2}-\lambda_1^{n+1/2})(\lambda_3^{n+1/2}-\lambda_2^{n+1/2}), \\
\delta^{n+1/2}\alpha_0^{n+1/2} & = e^{i\theta_1^{n+1/2}}\lambda_2^{n+1/2}\lambda_3^{n+1/2}(\lambda_3^{n+1/2}-\lambda_2^{n+1/2}) \\&+e^{i\theta_2^{n+1/2}}\lambda_3^{n+1/2}\lambda_1^{n+1/2}(\lambda_1^{n+1/2}-\lambda_3^{n+1/2}) \\
&+e^{i\theta_3^{n+1/2}}\lambda_1^{n+1/2}\lambda_2^{n+1/2}(\lambda_2^{n+1/2}-\lambda_1^{n+1/2}),\\
\delta^{n+1/2}\alpha_1^{n+1/2} & = i \Big(e^{i\theta_1^{n+1/2}}((\lambda_2^{n+1/2})^2-(\lambda_3^{n+1/2})^2) \\
& +e^{i\theta_2^{n+1/2}}((\lambda_3^{n+1/2})^2-(\lambda_1^{n+1/2})^2) \\
& +e^{i\theta_3^{n+1/2}}((\lambda_1^{n+1/2})^2-(\lambda_2^{n+1/2})^2)\Big),\\
\delta^{n+1/2}\alpha_2^{n+1/2} & = - \Big(e^{i\theta_1^{n+1/2}}(\lambda_3^{n+1/2}-\lambda_2^{n+1/2})\\
& +e^{i\theta_2^{n+1/2}}(\lambda_1^{n+1/2}-\lambda_3^{n+1/2}) \\& +e^{i\theta_3^{n+1/2}}(\lambda_2^{n+1/2}-\lambda_1^{n+1/2})\big),\\
\end{aligned}
\end{equation*}

Setting $\alpha^{n+1/2}=\alpha_0^{n+1/2}$, $\beta^{n+1/2}=-i\alpha_1^{n+1/2}$ and $\xi^{n+1/2} = \alpha_2^{n+1/2}/(\alpha_1^{n+1/2})^2$, the exact scheme for the Liouville equation becomes
\begin{equation}
\label{eq:exactLiouvillealt}
\begin{aligned}
\rho^{n+1} = &(\alpha^{n+1/2}I+\beta^{n+1/2}V^{n+1/2}+\xi^{n+1/2}(\beta^{n+1/2})^2(V^{n+1/2})^2)^{-1}\rho^n\\
&(\alpha^{n+1/2}I+\beta^{n+1/2}V^{n+1/2}+\xi^{n+1/2}(\beta^{n+1/2})^2(V^{n+1/2})^2)
\end{aligned}
\end{equation}
As we have seen for the two-level model in \cite{SBF17}, it is also possible to write the exact scheme for the three-level Liouville equation in the form of a NSFD model. Indeed the scheme \eqref{eq:exactLiouvillealt} can be written as
\begin{equation*}
\begin{aligned}
\alpha^{n+1/2}(\rho^{n+1}-\rho^n) 
= \beta^{n+1/2} & [\rho^n(V^{n+1/2}+\xi^{n+1/2}\beta^{n+1/2}(V^{n+1/2})^2)\\
- & (V^{n+1/2}+\xi^{n+1/2}\beta^{n+1/2}(V^{n+1/2})^2)\rho^{n+1}].
\end{aligned}
\end{equation*}
If $\dt$ is small enough, we can ensure that $\alpha^{n+1/2}$ and $\beta^{n+1/2}$ are nonzero and
\begin{equation*}
\begin{aligned}
i\dfrac{\alpha^{n+1/2}}{\beta^{n+1/2}}(\rho^{n+1}-\rho^n) 
= & -i[(V^{n+1/2}+\xi^{n+1/2}\beta^{n+1/2}(V^{n+1/2})^2) \rho^{n+1}\\
& -\rho^n (V^{n+1/2}+\xi^{n+1/2}\beta^{n+1/2}(V^{n+1/2})^2)]
\end{aligned}
\end{equation*}
or equivalently
\begin{equation}
\label{eq:schemainteraction3}
\left\{
\begin{array}{l}
\begin{aligned}
(\Phi^{n+1/2}(\dt))^{-1}(\rho^{n+1}-\rho^n) 
= & -i[(V^{n+1/2}+\xi^{n+1/2}\beta^{n+1/2}(V^{n+1/2})^2) \rho^{n+1}\\
& -\rho^n(V^{n+1/2}+\xi^{n+1/2}\beta^{n+1/2}(V^{n+1/2})^2)] \\
\Phi^{n+1/2}(\dt) = &-i\frac{\beta^{n+1/2}}{\alpha^{n+1/2}}I.
\end{aligned}
\end{array}
\right.
\end{equation}
In the left-hand side, we can recognize a nonstandard discretization in which the discretization time step-size undergoes a renormalization. And we can see that the renormalization matrix $\Phi$ has the following property:
\begin{equation}
\label{eq:property}
\Phi^{n+1/2}(\dt) = \dt I + \calO(\dt^2) \text{ when } \dt\rightarrow 0,
\end{equation}
because
\begin{equation*}
\label{eq:valimiteparam}
\lim_{\dt\rightarrow 0}\alpha^{n+1/2}=1;\
\lim_{\dt\rightarrow 0}\beta^{n+1/2}=i\dt;\
\lim_{\dt\rightarrow 0}\xi^{n+1/2} = \frac12.
\end{equation*}

The Strang splitting scheme derived from Equations \eqref{eq:SchemaRelaxNut} and \eqref{eq:exactLiouvillealt} is exactly \eqref{eq:canonicalscheme}. This scheme has a variable time step-size and preserve positiveness, because both steps \eqref{eq:SchemaRelaxNut} and \eqref{eq:exactLiouvillealt} are positive. The trace is also conserved.

\end{document}